# ABNORMAL SUBGROUPS AND CARTER SUBGROUPS IN SOME INFINITE GROUPS

L.A. Kurdachenko and I.Ya. Subbotin

**Abstract**. Some properties of abnormal subgroups in generalized soluble groups are considered. In particular, the transitivity of abnormality in metahypercentral groups is proven. Also it is proven that a subgroup *H* of a radical group G is abnormal in *G* if and only if every intermediate subgroup for *H* coincides with its normalizer in *G*. This result extends on radical groups the well-known criterion of abnormality for finite soluble groups due to D. Taunt. For some infinite groups (not only periodic) the existence of Carter subgroups and their conjugation have been also obtained.

**Keywords**: abnormal subgroups, pronormal subgroups, Carter subgroups.

**Mathematics Subject Classification**: 20E34, 20F19, 20F22

A subgroup $H$ of a group $G$ is abnormal in $G$ if $g \in \langle H, H^g \rangle$ for each element $g \in G$. Abnormal subgroups have appeared in the paper [HP] due to P. Hall, while the term "abnormal subgroup" itself belongs to R. Carter [CR]. Abnormal subgroups are antipodes to normal subgroups. In finite (especially soluble) groups the properties of abnormal subgroups have been studied in details. However we cannot say this about infinite groups. It is even unknown what groups contain proper abnormal subgroups. In this connection recall that a finite group is nilpotent if and only if it does not include proper abnormal subgroups. If $G$ is a locally nilpotent group, then $G$ has no proper abnormal subgroups [KUS]. However, in general, the converse assersion is not proven. In [KOS, KS 2] some classes of infinite groups, in which the absence of proper abnormal subgroups implies locally nilpotency, have been considered.

As normality, abnormality is not a transitive relation (see, for example group $S_4$). The groups (finite and infinite) with transitivity of normality are studied well enough (the most complete results one can find in [RD]). The situation for the groups with transitivity of abnormality is different. From the results of § 13 of Chapter VI of book [HB] it follows that in every finite metanilpotent group abnormality is transitive. For infinite groups this question has been considered in [KS 1, KS 2]. In the current article the following most general result has been obtained.

**1.2. Theorem**. *Let G be a group and suppose that A is a normal subgroup of G such that G/A has no proper abnormal subgroups. If A satisfies the normalizer condition, then in G abnormality is transitive.*

It is interesting to mention that this theorem and results of [RD] implies that soluble groups with transitivity of normality is a proper subclass of soluble groups in which abnormality is transitive.

In finite soluble groups abnormality is tightly connected to selfnormalizing. For example, D.Taunt has proved that a subgroup $H$ of a finite soluble group $G$ is abnormal if and only if every intermediate subgroup for $H$ coincides with its normalizer in $G$ (see, for example, [RD 6, 9.2.11]). Remind that *a subgroup S is said to be an intermediate subgroup for H if $H \leq S$*. The following theorem extends this result on radical groups.

**1.6. Theorem**. *Let G be a radical group and let H be a subgroup of G. Then H is abnormal in G if and only if every intermediate subgroup for H coincides with its normalizer in G.*

In finite groups there are many important families of subgroups having crucial influence on the group structure (for example Sylow and Hall subgroups, system normalizers, subgroups defined by formations, and so on). Many of them are very specific for the finite groups. With Carter subgroups the situation is different. In their original definition we can find no specifications related to finite groups. These subgroups are very tightly connected to abnormality. Indeed, Carter subgroups of finite groups can be defined as abnormal nilpotent subgroups. The following question naturally arises: which infinite groups posses Carter subgroups? S.E. Stonehewer in his papers [SE 1, SE 2] proved that periodic locally soluble groups having a locally nilpotent radical of finite index and locally soluble *FC* – groups are some examples of such infinite groups. In the current article we consider the following class of infinite groups.

Let $\mathfrak{X}$ is a class of groups. *A group G is said to be an artinian–by –$\mathfrak{X}$– group if G has a normal subgroup H such that G/H belongs to the class $\mathfrak{X}$ and H satisfies Min–G.*

This kind of groups has been introduced by D.J.S. Robinson [RD 3 – RD 5, RD 7, RD 8] and D.I. Zaitsev [ZD 1 – ZD 5] in their series of papers dedicated to the existence of complements to the $\mathfrak{X}$ – residual (when it is abelian) for some natural classes $\mathfrak{X}$ (such as hypercentral groups, locally nilpotent groups, hypercyclic groups, locally supersoluble groups, hyperfinite groups). For the classes $\mathfrak{X}$ considered in our paper this definition implies that if $G$ is an artinian – by –$\mathfrak{X}$– group and $R = G_{\mathfrak{X}}$ is its $\mathfrak{X}$– residual, then $G/R \in \mathfrak{X}$ and $R$ satisfies Min–$G$.

We will deal with artinian – by – hypercentral groups whose locally nilpotent residual is nilpotent. It is a natural first step. Since these groups are generalizations of finite metanilpotent groups, we will use for the definition some characterizations of Carter subgroups, which are valid for finite metanilpotent groups. In particular, in finite metanilpotent groups Carter subgroups coincide with system normalizers (see, for example, [RD 6, 9.5.10]). By a result due to P. Hall (see, for example, [RD 6, 9.2.15]) the system normalizers of a finite soluble group are precisely its minimal subabnormal subgroups. As we have noted above, in a finite metanilpotent group every subabnormal subgroup is abnormal. Consequently, we can define a Carter subgroup in a metanilpotent group as a minimal abnormal subgroup.

**2.1. Theorem**. *Let G be an artinian – by – hypercentral group and suppose that its locally nilpotent residual K is nilpotent. Then G includes a minimal abnormal subgroup L. Moreover, L is a maximal hypercentral subgroup and it includes the upper hypercenter of G. In particular, G = KL. If H is another minimal abnormal subgroup, then H is conjugate to L.*

**2.2. Corollary**. *Let G be an artinian – by – hypercentral group and suppose that its locally nilpotent residual K is nilpotent. Then G includes a hypercentral abnormal subgroup L. Moreover, L is maximal hypercentral subgroup and it includes the upper hypercenter of G. In particular, G = KL. If H is another hypercentral abnormal subgroup, then H is conjugate to L.*

Let *G* be an artinian – by – hypercentral group with a nilpotent hypercentral residual.

*A subgroup L is called a Carter subgroup of a group G if H is a hypercentral abnormal subgroup of G (or, equivalently, H is a minimal abnormal subgroup of G).*

A Carter subgroup in finite soluble group can be defined as a covering subgroup for the formation of nilpotent groups. As we will see, this characterization can be extended on the groups under consideration.

Recall that a subgroup *H* of a group *G* is said to be *a $\mathbf{L}\mathfrak{N}$ -covering subgroup if H is locally nilpotent and if $S = S_{\mathbf{L}\mathfrak{N}}H$ for every subgroup S, which includes H. Here $S_{\mathbf{L}\mathfrak{N}}$ is a locally nilpotent residual of subgroup S.*

**2.3. Theorem**. *Let G be an artinian – by – hypercentral group and suppose that its locally nilpotent residual K is nilpotent. If L is a Carter subgroup of G, then L is a $\mathbf{L}\mathfrak{N}$ – covering subgroup of G. Conversely, if H is a $\mathbf{L}\mathfrak{N}$ – covering subgroup of G, then H is a Carter subgroup of G.*

In a finite soluble group $\mathfrak{N}$– covering subgroups are exactly $\mathfrak{N}$ – projectors. Therefore a Carter subgroup of a finite soluble group can be defined as an $\mathfrak{N}$ – projector. As we will see, this characterization also can be extended on the considered in this article groups.

Recall that *a subgroup L of a group G is said to be a locally nilpotent projector, if LH/H is a maximal locally nilpotent subgroup of G/H for each normal subgroup H of a group G.*

**2.4. Theorem**. *Let G be an artinian – by – hypercentral group and suppose that its locally nilpotent residual K is nilpotent. If L is a Carter subgroup of G, then L is a locally nilpotent projector of G. Conversely, if D is a locally nilpotent projector of G, then H is a Carter subgroup of G.*

# 1. Abnormal subgroups in some infinite groups.

We will consider some classes of groups, in which abnormality is transitive.

**1.1. Lemma**. *Let G be a group and B be an abnormal subgroup of G. Suppose that B = HA where A is a normal subgroup of a group G. If H is an abnormal subgroup of B, then H is an abnormal subgroup of G.*

This assertion belongs to P. Hall and one can find it, for example, in the book [RD 6, 9.2.12] where it is proven for finite groups $G$. However this proof does not use the finiteness of G and the result is valid for infinite groups also.

**1.2. Theorem**. *Let G be a group and suppose that A is a normal subgroup of G such that G/A has no proper abnormal subgroups. If A satisfies the normalizer condition, then in G abnormality is transitive.*

**Proof**. Let $B$ be an abnormal subgroup of $G$ and let $H$ be an abnormal subgroup of $B$. If $B = G$, then all is proved. Suppose that $B \neq G$. By elementary properties of abnormal subgroups $AB$ is abnormal in $G$. Since $G/A$ has no proper abnormal subgroup, $AB = G$. Proceeding by induction, we will construct a strictly ascending chain

$$B = B_0 < B_1 < \ldots B_\alpha < B_{\alpha+1} < \ldots$$

of subgroups such that $H$ is abnormal in every subgroup $B_\alpha$. Put $K = A \cap B$. Then $K$ is a normal subgroup of $B$, so that $KH$ is abnormal in $B$ and $KH/K$ is abnormal in $B/K$. Since $G/A = AB/A \cong B/(A \cap B) = B/K$, we obtain that $B/K$ has no proper abnormal subgroups. It follows that $B/K = HK/K$, that is $B = HK$. Since $A$ satisfies the normalizer condition, there exists a subgroup $K_1$ such that $K < K_1$ and $K$ is a normal subgroup of $K_1$. Put $B_1 = \langle B, K_1 \rangle$. Then $B \neq B_1$. Since $K$ is normal in both $B$ and $K_1$, $K$ is normal in $B_1$. By Lemma 1.1 the equation $B = HK$ implies that $H$ is abnormal in $B_1$. Suppose that we have already constructed the subgroups $B_\beta$ for all $\beta < a$. Let first $\alpha - 1$ exists and $B_{\alpha-1} \neq G$. Since $B \leq B_{\alpha-1}$, the latter is abnormal in $G$. By elementary properties of abnormal subgroups $AB_{\alpha-1}$ is abnormal in $G$. Since $G/A$ has no proper abnormal subgroups, $AB_{\alpha-1} = G$. Put $K_{\alpha-1} = A \cap B_{\alpha-1}$. Then $K_{\alpha-1}$ is a normal subgroup of $B_{\alpha-1}$, so that $K_{\alpha-1}H$ is abnormal in $B_{\alpha-1}$ and $K_{\alpha-1}H/K_{\alpha-1}$ is abnormal in $B_{\alpha-1}/K_{\alpha-1}$. However, since $G/A = AB_{\alpha-1}/A \cong B_{\alpha-1}/(A \cap B_{\alpha-1}) = B_{\alpha-1}/K_{\alpha-1}$, we obtain that $B_{\alpha-1}/K_{\alpha-1}$ has no proper abnormal subgroups. It follows that $B_{\alpha-1}/K_{\alpha-1} = K_{\alpha-1}H/K_{\alpha-1}$, that is $B_{\alpha-1} = K_{\alpha-1}H$. Since $A$ satisfies the normalizer condition, there exists a subgroup $K_\alpha$ such that $K_{\alpha-1} < K_\alpha$ and $K_{\alpha-1}$ is a normal subgroup of $K_\alpha$. Put $B_\alpha = \langle B_{\alpha-1}, K_\alpha \rangle$. Then $B_{\alpha-1} \neq B_\alpha$. Since $K_{\alpha-1}$ is normal in both $B_{\alpha-1}$ and $K_\alpha$, $K_{\alpha-1}$ a is normal in $B_\alpha$. By Lemma 1.1 equation $B_{\alpha-1} = K_{\alpha-1}H$ implies that $H$ is abnormal in $B_\alpha$.

Suppose now that $\alpha$ is a limit ordinal. Then put $B_\alpha = \bigcup_{\beta < \alpha} B_\beta$. Let $x$ be an arbitrary element of $B_\alpha$. Then $x \in B_\nu$ for some $\nu < \alpha$. We have $\langle H, H^x \rangle \leq B_\nu$. Since $H$ is abnormal in $B_\nu$, $x \in \langle H, H^x \rangle$. This means that $H$ is an abnormal subgroup of $B_\alpha$.

By our construction, $B_{\alpha-1} \neq B_\alpha$ for each ordinal $\alpha$. It follows that there is an ordinal $\gamma$ such that $B_\gamma = G$. Thus $H$ is abnormal in $B_\gamma = G$. □

**1.3. Corollary**. *Let G be a group and suppose that A is a normal subgroup of G such that G/A does not include proper abnormal subgroup. If A is hypercentral, then in G abnormality is transitive.*

**1.4. Corollary**. *Let G be a group and suppose that A is a normal subgroup of G such that G/A is locally nilpotent. If A is hypercentral, then in G abnormality is transitive.*

Indeed, a locally nilpotent group does not include proper abnormal subgroups [KUS].

**1.5. Lemma**. *Let G be a group, H be a subgroup of G and let D be an H – invariant subgroup. Suppose that every intermediate subgroup for H coincides with its normalizer in G. If D is a locally nilpotent subgroup, then H is abnormal in HD.*

**Proof**. Put $L = HD$. Choose an arbitrary element $d \in D$ and consider the subgroup $K = \langle H, H^d \rangle$. Since $H \leq K, L = DK$. It follows that $B = D \cap K$ is a normal subgroup of $K$, in particular, $B$ is $H$ – invariant. If $d \in B$, then $d \in HB = K = \langle H, H^d \rangle$. Suppose that $d \notin B$. In the subgroup $U = \langle d, B \rangle$ choose a subgroup $M$, which is maximal with the properties $B \leq M$ and $d \notin M$. Clearly $M$ is a maximal subgroup of $U$. Since $U$ is locally nilpotent, $M$ is normal in $U$ ( see, for example, [RD 1, Theorem 5.38]). In particular, $d^{-1}Md = M$. If $h \in H$ then $[h, d] \in K$, that is $[h, d] \in D \cap K$. Since $H \leq K, K = (D \cap K)H = BH$. We have now $h^d = h[h, d] \in HB = K$ for each element $h \in H$. Let $y \in M, h \in H$ and consider the element $d^{-1}(h^{-1}yh)d = d^{-1}h^{-1}yhd$. Since $dhd^{-1}h^{-1} = b \in B$, $d^{-1}h^{-1} = h^{-1}d^{-1}b$. Similarly, $hd = adh$ for some element $a \in B$. Now we have $d^{-1}h^{-1}yhd = h^{-1}d^{-1}(bya)dh$. The inclusion $B \leq M$ implies that $bay \in M$. Then $d^{-1}(bya)d \in M$ and hence $h^{-1}d^{-1}(bya)dh. \in h^{-1}Mh$. In other words, $d^{-1}(h^{-1}Mh)d = h^{-1}Mh$. Let $C = \bigcap_{h \in H}(h^{-1}Mh.)$, then by above $d^{-1}Cd = C$. Since $B \leq M, B = h^{-1}Bh. \leq h^{-1}Mh$ for each $h \in H$, so that $B \leq C$. Furthermore, $d^{-1}Hd \leq K = HB \leq HC$.

It follows that $d^{-1}(HC)d \leq HC$. In other words, $d \in N_G(HC)$. Since $HC \cap D = C(H \cap D) \leq CB = C, d \notin HC$. On the other hand, by the conditions of our lemma $HC$ is self – normalizing. This contradiction shows that $d \in B$ and hence $d \notin K = \langle H, H^d \rangle$ what means that $H$ is abnormal in $HD$. □

**1.6. Theorem**. *Let G be a radical group and let H be a subgroup of G. Then H is abnormal in G if and only if every intermediate subgroup for H coincides with its normalizer in G.*

**Proof**. Let

$$\langle 1 \rangle = A_0 \leq A_1 \leq \ldots A_\alpha \leq A_{\alpha+1} \leq \ldots A_\gamma = G$$

be a series of normal subgroups of $G$, whose factors are locally nilpotent. If $H$ is abnormal in $G$, then every intermediate subgroup for $H$ is self – normalizing (see, for example, [RD 6, p. 266]). Conversely, suppose that every intermediate subgroup for $H$ is self – normalizing and will prove that $H$ is abnormal in $G$. More precisely, we will prove that $H$ is abnormal in $A_\alpha H$ for each $\alpha \leq \gamma$. By Lemma 1.5 $H$ is abnormal in $HA_1$, so, for the case $\alpha = 1$ all is proved. Suppose that we have already proved that $H$ is abnormal in $A_\beta H$ for all $\beta < \alpha$. Choose an arbitrary element $a \in A_\alpha$ and consider the subgroup $K = \langle H, H^a \rangle$. Suppose firstly that $\alpha$ is a limit ordinal. Then there is an ordinal $\beta < \alpha$ such that $a \in A_\beta$. Since $H$ is abnormal in $A_\beta H$, $a \in K = \langle H, H^a \rangle$. It follows that $H$ is abnormal in $A_\alpha H$. Assume now that $\alpha$ is not limit. Put $U = A_\alpha$ and $V = A_{\alpha-1}$. Consider the factor – group $UH/V$. Let $Z/V$ be an intermediate subgroup for $HV/V$ in $HU/V$. If $xV$ is an element of $UH/U$ such that $Z/V = (Z/V)^{xV}$, then $(Z/V)^{xV} = (Z^xV)/V = Z^x/V$ implies $Z^x = Z$. Since $VH \leq Z$, the subgroup $Z$ is self – normalizing, in particular, $x \in Z$. In other words, every subgroup of $UH/V$, which is intermediate for $HV/V$, is self – normalizing. By Lemma 1.5 $HV/V$ is abnormal in $HU/V$. In turn it follows that $HV$ is abnormal in $HU$. By the induction hypothesis $H$ is abnormal in $VH$. The application of Lemma 1.1 gives that $V$ is abnormal in $HU = HA_\alpha$. Since it is valid for each ordinal $\alpha$, $V$ is abnormal in $HA_\gamma = G$. □

**1.7. Corollary**. *Let G be a hyperabelian group and let H be a subgroup of G. Then H is abnormal in G if and only if every intermediate subgroup for H coincides with its normalizer in G.*

**1.8. Corollary**. *Let G be a soluble group and let H be a subgroup of G. Then H is abnormal in G if and only if every intermediate subgroup for H coincides with its normalizer in G.*

For the consideration of artinian – by – $\mathfrak{X}$ – groups we need some notions connected with the acting of a group on its abelian subgroup. It is better to formulate these concepts using the Modules Theory language.

Let $R$ be a ring, $G$ be a group, $A$ be an $RG$- module, $B, C$ be $RG$ - submodules of $A$ such that $B \leq C$. Factor $C/B$ is called central (more exactly $RG$ – central) (respectively eccentric or $RG$ – eccentric) if $G = C_G(C/B)$ (respectively $G \neq C_G(C/B)$). Put

$$\zeta_{RG}(A) = \{a \in A | aRG \text{ is } RG\text{ –central}\}.$$

Clearly $\zeta_{RG}(A)$ is an $RG$ – submodule of $A$. This submodule is called *the RG – center of A*.

Starting from the $RG$ - center we can construct *the upper RG – central series of the module A*:

$$\langle 0 \rangle = A_0 \leq A_1 \leq \ldots A_\alpha \leq A_{\alpha+1} \leq \ldots A_\gamma,$$

where $A_1 = \zeta_{RG}(A), A_{\alpha+1}/A_\alpha = \zeta_{RG}(A/A_\alpha), \alpha < \gamma$, and $\zeta_{RG}(A/A_\gamma) = \langle 0 \rangle$.

*The last term $A_\gamma$ of this series is called the upper RG – hypercenter of the module A and denote by $\zeta_{RG}^\infty(A)$. If $A = A_\gamma$, then the module A is called RG – hypercentral; if $\gamma$ is finite, then A is called RG – nilpotent.*

On the other hand, an *RG – submodule C of A is said to be RG - hypereccentric, if it has an ascending series*

$$\langle 0 \rangle = C_0 \leq C_1 \leq \ldots C_\alpha \leq C_{\alpha+1} \leq \ldots C_\gamma = C$$

*of RG– submodules of A such that each factor $C_{\alpha+1}/C_\alpha$ is an RG – eccentric simple RG – module, for every $\alpha < \gamma$.*

Following D.I. Zaitsev [ZD 1], we say that *an RG-module A has the Z - decomposition if the following equality holds $A = \zeta_{RG}^\infty(A) \oplus \varepsilon_{RG}^\infty(A)$ where $\varepsilon_{RG}^\infty(A)$ is the maximal RG - hypereccentric RG - submodule of A.*

Note that if $A$ is an artinian $RG$ – module, then $\varepsilon_{RG}^\infty(A)$ includes every $RG$ – hypereccentric $RG$ – submodule, in particular, it is unique. In fact, let $B$ be an $RG$ – hypereccentric $RG$ – submodule of $A, E = \varepsilon_{RG}^\infty(A)$. If $(B + E)/E$ is non – zero it includes a non – zero simple $RG$ – submodule $U/E$. Since $(B + E)/E \cong B/(B \cap E), U/E$ is $RG$ – isomorphic with some simple $RG$ – factor of $B$ and it follows that $G \neq C_G(U/E)$. On the other hand, $(B + E)/E \leq A/E \cong \zeta_{RG}^\infty(A)$, that is $G = C_G(U/E)$. This contradiction proves that $B \leq E$. Hence $\varepsilon_{RG}^\infty(A)$ includes every $RG$ - hypereccentric $RG$ – submodule, in particular, it is unique.

**1.9. Lemma**. *Let G be a group and suppose that A is a normal abelian subgroup of G such that G/A is hypercentral and A satisfies Min −G. If $\varepsilon_{RG}^\infty(A) \neq \langle 1 \rangle$ and H is a complement to $\varepsilon_{RG}^\infty(A)$, then H is abnormal in G.*

**Proof**. By Theorem 1´ of [ZD 1] $A$ has the $Z$ – decomposition $A = C \times E$, where $C = \zeta_{ZG}^\infty(A), E = \varepsilon_{RG}^\infty(A)$. Factor – group $G/E$ is hypercentral, so that by Theorem 2 of [ZD 2] $G$ includes a subgroup $H$ such that $G = EH$ and $E \cap H = \langle 1 \rangle$. Moreover, every complement to $E$ in $G$ is conjugate to $H$. Let $S$ be an intermediate subgroup for $H$, that is $H \leq S$. Then $S = (E \cap S)H$. Clearly, $E \cap S$ is a $G$ – invariant subgroup, moreover, every $S$ – invariant (even $H$ – invariant) subgroup of $E \cap S$ is $G$ – invariant. Put $L = N_G(S)$. Similarly $L = (E \cap L)H$. If $L \neq S$, then $E \cap S \neq E \cap L$. In factor – group $L/(E \cap S)$ we have $S/(E \cap S) = H(E \cap S)/(E \cap S)$. If $x \in (L \cap E) \backslash S$, then $[x, S] \leq S$, so that $[x, S](E \cap S)/(E \cap S) \leq H(E \cap S)/(E \cap S)$. On the other hand, $[x, S] \leq E$, that is $[x, S] \leq E \cap S$. In other words, factor $(E \cap L)/(E \cap S)$ is $S$ – central, and therefore $G$ – central. However, by the selection of $E$ every $G$ – chief factor of $E$ is $G$ – eccentric, so that $E$ does not have the $G$ - central factors. This contradiction shows that

$S = N_G(S)$. By Corollary 1.7 $H$ is an abnormal subgroup of $G$. □

The next theorem is about artinian – by – hypercentral groups with no proper abnormal subgroups.

**1.10. Theorem**. *Let G be a group and suppose that T is a normal soluble subgroup of G such that G/T is hypercentral and T satisfies Min –G. Group G is hypercentral if and only if G has no proper abnormal subgroups.*

**Proof**. Let

$$\langle 1 \rangle = T_0 \leq T_1 \leq \ldots \leq T_d = T$$

be the derived series of $T$. We will use induction by $d$. Let $d = 1$, that is $A = T_1$ is abelian. Then $A$ has the $Z$– decomposition $A = C \times E$, where $C = \zeta_{RG}^\infty(A)$, $E = \varepsilon_{RG}^\infty(A)$ [ZD 1, Theorem 1´]. Suppose that $G$ is not hypercentral. It follows that $E \neq \langle 1 \rangle$. Factor – group $G/E$ is hypercentral, so that by Theorem 2 of [ZD 2] $G$ includes a subgroup $H$ such that $G = EH$ and $E \cap H = \langle 1 \rangle$. By Lemma 1.14 subgroup $H$ is abnormal in $G$. This contradiction shows that $E = \langle 1 \rangle$, and hence $G$ is hypercentral.

Let now $d > 1$ and we have already proved that $G/T_1$ is hypercentral. We can repeat the above arguments and obtain that $G$ is hypercentral. □

## 2. Carter subgroups in some infinite groups.

We will consider the existence of Carter subgroups in some artinian – by – hypercentral groups.

Note that if $G$ is an artinian – by – hypercentral group, then its hypercentral residual $R$ coincides with locally nilpotent residual $K$, moreover $G/K$ is hypercentral. Indeed, $K \leq R$. Since $G/K$ is locally nilpotent, it has a central series $\mathfrak{Z}$ (see, for example, [RD 2, Corollary to Theorem 8.24]). Put $Z_1 = \{Z \cap (R/K) | Z \in \mathfrak{Z}\}$. Since $R/K$ satisfies Min – G, it has a minimal element $C/K$. But every chief factor of a locally nilpotent group is central (see, for example, [RD 1, Corollary 1 to Theorem 5.27]). In other words, $\zeta(G/K)$ is non – identity. Using the similar arguments and transfinite induction, we obtain that $R/K$ has an ascending $G$ – central series. Since $G/R$ is hypercentral, $G/K$ is likewise hypercentral.

**2.1. Theorem**. *Let G be an artinian – by – hypercentral group and suppose that its locally nilpotent residual K is nilpotent. Then G includes a minimal abnormal subgroup L. Moreover, L is a maximal hypercentral subgroup and it includes the upper hypercenter of G. In particular, $G = KL$. If H is another minimal abnormal subgroup, then H is conjugate to L.*

**Proof**. Let

$$\langle 1 \rangle = K_0 \leq K_1 \leq \ldots \leq K_c = K$$

be the upper central series of $K$. We will use induction by $c$. Let $c = 1$, that is $A = K_1$ is abelian. Then $A$ has the $Z$ – decomposition $A = C \times E$, where $C = \zeta_{RG}^\infty(A)$, $E = \varepsilon_{RG}^\infty(A)$ [ZD 1, Theorem 1´]. Factor – group $G/E$ is hypercentral, so that by Theorem 2 of [ZD 2] $G$ includes a subgroup $L$ such that $G = EL$ and $E \cap L = \langle 1 \rangle$. By isomorphism $L \cong G/E$, subgroup $L$ is hypercentral. By Lemma 1.9 $L$ is abnormal in $G$. Let $S$ be an intermediate subgroup for $L$. Then $S = (E \cap S)L$. Since $E$ is abelian, the equation $G = ES$ implies that every $S$ – invariant subgroup of $E$ is also $G$ – invariant. In particular, $E \cap S$ satisfies Min – $G$. Every $G$ – chief factor of $E$ is $G$ – eccentric. Thus every $S$ – chief factor of $E \cap S$ is $S$ – eccentric. It follows that if $S \neq L$ (that is $E \cap S \neq \langle 1 \rangle$), then $S$ can not be a hypercentral subgroup. Let $H$ be another minimal abnormal subgroup of $G$. Then $HE/E$ is an abnormal subgroup of $G/E$. Since a hypercentral group does not include proper abnormal subgroups [KUS], $HE/E = G/E$ or $HE = G$. Suppose that $E \cap H \neq \langle 1 \rangle$. Again each $H$ - invariant subgroup of $E$ is also $G$ – invariant. In particular, $E \cap H$ satisfies Min –$H$ and $E \cap H = \varepsilon_{RG}^\infty(E \cap H)$. By Theorem 2 of [ZD 2] $H$ includes a subgroup $U$ such that $H = (E \cap H)U$ and $(E \cap H) \cap U = \langle 1 \rangle$. Since $E \cap H \neq \langle 1 \rangle$, $H \neq U$. By Lemma 1.9 $U$ is abnormal in $H$. Theorem 1.2 yields that $U$ is abnormal in $G$. However, this contradicts the selection of $H$. This contradiction shows that $E \cap H = \langle 1 \rangle$. By Theorem 2 of [ZD 2] $H$ is conjugate to subgroup $L$.

Suppose now that $c > 1$ and we have already proved theorem for factor – group $G/A$. Let $V/A$ be a minimal abnormal subgroup of $G/A$. Since $C_G(A) \geq T$, $G/C_G(A)$ is hypercentral. Then $A$ has the $Z$ – decomposition $A = C \times E$, where $C = \zeta_{RG}^\infty(A)$, $E = \varepsilon_{RG}^\infty(A)$ [ZD 1, Theorem 1´]. Equation $G/A = (T/A)(V/A)$ and inclusion $A \leq \zeta(T)$ imply that every $V$ – invariant subgroup of $A$ is $G$ – invariant. In other words, $A$ satisfies Min – $V$ and $C = \zeta_{RG}^\infty(A)$, $E = \varepsilon_{RG}^\infty(A)$. Factor – group $V/E$ is hypercentral, so that by Theorem 2 of [ZD 2] $V$ includes a subgroup $L$ such that $V = EL$ and $E \cap L = \langle 1 \rangle$. By isomorphism $L \cong V/E$, $L$ is hypercentral. As we proved above, $L$ is abnormal in $V$, and Theorem 1.2 yields that $L$ is abnormal in $G$ because $V$ is abnormal in $G$ by induction hypothesis. Since $L$ is hypercentral, $L$ is a minimal abnormal subgroup of $G$. Let $H$ be another minimal abnormal subgroup of $G$. Since $G/T$ is hypercentral, $HT = G$. Inclusion $T_1 \geq \zeta(T)$ implies that every $H$ – invariant subgroup of $T_1$ is also $G$ – invariant. In particular, $T_1$ satisfies Min – $H$. By the same reason, every factor $T_j/T_{i-1}$ satisfies Min – $H$ for each $j, 1 \leq j \leq c$. It follows that $T$ satisfies Min –$H$, in particular, $H \cap T$ satisfies Min – $H$. By isomorphism $H/(H \cap T) \cong HT/T$ factor – group $H/(H \cap T)$ is hypercentral. If $W$ is an abnormal subgroup of $H$, then by Theorem 1.2 $W$ is abnormal in $G$. By the selection of $H$ it follows that $H$ does not include a proper abnormal (in $H$) subgroup. By Theorem 1.10 $H$ is hypercentral. So $HA/A$ is an abnormal hypercentral subgroup of $G/A$. On other words, $HA/A$ is a minimal abnormal subgroup of $G/A$. By inductive hypothesis there is an element $g \in G$ such that $(LA/A)^{gA} = HA/A$. Then $H \leq L^{gA}$. Since every $H$ – invariant subgroup of $A$ is also $G$ – invariant, $A$ satisfies Min – $H$ and $C = \zeta_{RG}^\infty(A)$, $E = \varepsilon_{RG}^\infty(A)$. In particular, $HC$ is hypercentral. Since $H$ is abnormal in $HC$, it follows that $HC = H$. Furthermore, $H \cap E = \langle 1 \rangle$, so $H$ is a complement to $E$ in $HA = L^{gA}$. By Theorem 2 of [ZD 2] there is an element $z \in HA$ such that

$H = L^{gz}$.   □

**2.2. Corollary**. *Let G be an artinian – by – hypercentral group and suppose that its locally nilpotent residual K is nilpotent. Then G includes a hypercentral abnormal subgroup L. Moreover, L is maximal hypercentral subgroup and it includes the upper hypercenter of G. In particular, G = KL. If H is another hypercentral abnormal subgroup, then H is conjugate to L.*

By Theorem 2.1 a minimal abnormal subgroup is hypercentral. Conversely, if $L$ is an abnormal hypercentral subgroup then $L$ does not include a proper abnormal subgroup [KUS]. This means that $L$ is a minimal abnormal subgroup.

**2.3. Theorem**. *Let G be an artinian – by – hypercentral group and suppose that its locally nilpotent residual K is nilpotent. If L is a Carter subgroup of G, then L is a $\mathbf{L}\mathfrak{N}$ – covering subgroup of G. Conversely, if H is a $\mathbf{L}\mathfrak{N}$ – covering subgroup of G, then H is a Carter subgroup of G.*

**Proof**. Let $L$ be a Carter subgroup of $G$ (its existence follows from Theorem 2.1). If $L \leq S$, then $L$ is an abnormal subgroup of $S$. Let

$$\langle 1 \rangle = K_0 \leq K_1 \leq \ldots \leq K_c = K$$

be the upper central series of $T$. By Theorem 2.1 $G = LK$. It follows that every $L$ – invariant subgroup of $K_j/K_{j-1}$ is $G$ – invariant, $1 \leq j \leq c$. In particular, $K_j/K_{j-1}$ satisfies Min – $L$, $1 \leq j \leq c$. Put $Q_j = S \cap K_j, 0 \leq j \leq c$. Then every factor $Q_j/Q_{j-1}$ is central in $S \cap K$ and satisfies Min – $L$. It follows that $S \cap K$ is a nilpotent subgroup satisfying Min – $L$. Since $S/(S \cap K) \cong SK/K$ is hypercentral and $L(S \cap K)/(S \cap K)$ is its abnormal subgroup, $L(S \cap K) = S$, because a hypercentral group does not include a proper abnormal subgroup. It follows that $S \cap K$ satisfies Min –$S$. In other words, $S$ is an artinian – by – hypercentral group. Let $R$ be the locally nilpotent residual of $S$. As we have already noted prior to this theorem, $S/R$ is hypercentral. Since $LR/R$ is an abnormal subgroup of $S/R$, $LR = S$. This means that $L$ is a $\mathbf{L}\mathfrak{N}$ – covering subgroup of $G$.

Conversely, let $H$ be an arbitrary $\mathbf{L}\mathfrak{N}$ – covering subgroup of $G$. We will prove that $H$ is an abnormal hypercentral subgroup of $G$. We will use for this induction by $c$. Let first $c = 1$, that is $A = K_1$ is abelian. Then $A$ has the $Z$ – decomposition $A = C \times E$, where $C = \zeta_{RG}^\infty(A)$, $E = \varepsilon_{RG}^\infty(A)$ [ZD 1, Theorem 1´]. Clearly, $E$ is the hypercentral residual (and the locally nilpotent residual) of $G$, so that $EH = G$. Suppose that $B = E \cap H \neq \langle 1 \rangle$. Obviously, $B$ is a $G$ – invariant subgroup of $E$. Since $E$ satisfies Min – $G$, $B$ includes a minimal $G$ – invariant subgroup $M$. Equation $G = EH$ yields that every $H$ – invariant subgroup of $E$ is likewise $G$ – invariant. It follows that $M$ is a minimal $H$ – invariant subgroup. However every chief factor of a locally nilpotent group $H$ is central ( see, for example, [RD 1, Corollary 1 to Theorem 5.27]), so that $M$ is $H$ - central. In this case $M$ is central in $G$. This contradicts the inclusion $M \cap E = \varepsilon_{RG}^\infty(A)$. This contradiction shows that $H \cap E = \langle 1 \rangle$. Then $H \cong HE/H = G/E$ is hypercentral. By Lemma 1.9 *H* is abnormal in *G*.

Suppose now that $c > 1$ and consider factor – group *G/A*. It is not hard to see that *HA/A* is a **L𝔑** – covering subgroup of *G/A*. By induction hypothesis *HA/A* is an abnormal hypercentral subgroup. Consider now a subgroup *HA*. Clearly *H* is a **L𝔑** – covering subgroup of *HA*. We have already proved that *H* is a hypercentral subgroup and *H* is abnormal in *HA*. By Lemma 1.1 *H* is abnormal in *G*. Hence *H* is a Carter subgroup of *G*. □

**2.4. Theorem**. *Let G be an artinian – by – hypercentral group and suppose that its locally nilpotent residual K is nilpotent. If L is a Carter subgroup of G, then L is a locally nilpotent projector of G. Conversely, if D is a locally nilpotent projector of G, then H is a Carter subgroup of G.*

**Proof**. Let *L* be a Carter subgroup of *G*; its existence follows from Theorem 2.1. If *H* is a normal subgroup of *G*, then *LH/H* is an abnormal subgroup of *G/H*. Let *K/H* be a locally nilpotent subgroup including *LH/H*. Since a locally nilpotent group does not include a proper abnormal subgroups [KUS], *LH/H* = *K/H*. This means that *L* is a locally nilpotent projector.

Conversely, let *D* be an arbitrary locally nilpotent projector of *G*. Since *DK/K* is a maximal locally nilpotent subgroup of a hypercentral group *G/E*, *DK = G*. Let

$$\langle 1 \rangle = K_0 \leq K_1 \leq \ldots \leq K_c = K$$

be the upper central series of *T*. We will prove that *D* is an abnormal hypercentral subgroup of *G*. We will use for this the induction by *c*. Let first $c = 1$, that is $A = K_1$ is abelian. Then *A* has the *Z* – decomposition $A = C \times E$, where $C = \zeta_{RG}^\infty(A)$, $E = \varepsilon_{RG}^\infty(A)$ [ZD 1, Theorem 1´]. Since *E* is a hypercentral residual of *G*, by above *DE = G*. As in the proof of Theorem 2.3 we can prove that $E \cap D = \langle 1 \rangle$. Then $D \cong DE/H = G/E$ is hypercentral. By Lemma 1.9 *D* is abnormal in *G*. Suppose now that $c > 1$ and consider factor – group *G/A*. Obviously *DA/A* is a locally nilpotent projector of *G/A*. By induction hypothesis *DA/A* is an abnormal hypercentral subgroup. Consider now a subgroup *HA*. The inclusion $A \leq \zeta(K)$ and equation $G = KD$ implies that every *D* – invariant subgroup of *A* is *G* – invariant. In particular, *A* satisfies Min – *G*. Since *DA/A* is hypercentral, *A* has the *Z* – decomposition $A = Z \times U$, where $C = \zeta_{RD}^\infty(A)$, $U = \varepsilon_{RD}^\infty(A)$ [ZD 1, Theorem 1']. Moreover, $C = \zeta_{RG}^\infty(A) = \zeta_{RD}^\infty(A)$, $\varepsilon_{RG}^\infty(A) = \varepsilon_{RD}^\infty(A)$. Since *D* is a maximal locally nilpotent subgroup, $Z \leq D$. Using the arguments of the proof of Theorem 2.3 and the equation $\varepsilon_{RG}^\infty(A) = \varepsilon_{RD}^\infty(A)$ we can prove that $D \cap U = \langle 1 \rangle$. Lemma 1.9 implies that *D* is abnormal in *DA*. By Lemma 1.1 *D* is abnormal in *G*. Hence *H* is a Carter subgroup of *G*. □